\documentclass[11pt]{article}
\usepackage{amsmath,amssymb,amsthm,geometry,verbatim,enumerate,float,tikz}
\newtheorem{thm}{Theorem}

\newtheorem{coro}[thm]{Corollary}

\usepackage{setspace}

\begin{document}

\onehalfspace

\title{Independence and Matching Number in Graphs with Maximum Degree 4}

\author{Felix Joos}

\date{}

\maketitle

\vspace{-1cm}

\begin{center}
Institut f\"{u}r Optimierung und Operations Research, 
Universit\"{a}t Ulm, Ulm, Germany\\
felix.joos@uni-ulm.de
\end{center}

\begin{abstract}
We prove that $\frac{7}{4}\alpha(G)+\beta(G)\geq n(G)$ and $\alpha(G)+\frac{3}{2}\beta(G)\geq n(G)$ for every triangle-free graph $G$ with maximum degree at most $4$,
where $\alpha(G)$ is the independence number and $\beta(G)$ is the matching number of $G$, respectively.
These results are sharp for a graph on $13$ vertices.
Furthermore we show $\chi(G)\leq \frac{7}{4}\omega(G)$ for $\{3K_1,K_1\cup K_5\}$-free graphs, 
where $\chi(G)$ is the chromatic number and $\omega(G)$ is the clique number of $G$, respectively.
\end{abstract}



\maketitle

\section{Introduction}
\noindent
Inspired by a result of Henning, L\"owenstein, and Rautenbach \cite{helora} 
we investigate  linear combinations of the independence number and the matching number of graphs.

The intuition is that the matching number and the independence number are negatively correlated,
that is, if the matching number is small, then the independence number is large and vice versa.
Thus the sum of both is bounded from below.

We start with some notation.
We only consider simple, finite, and undirected graphs.
For a graph $G$, we denote by $V(G)$ and $E(G)$ the vertex and edge set, respectively.
Let $n(G)=|V(G)|$ be the order and $m(G)=|E(G)|$ be the size of $G$.
For a vertex $v\in V(G)$, let $N_G(v)$ be the set of all neighbours of $v$.
Define the degree $d_G(v)$ of $v$ by $d_G(v)=|N_G(v)|$.
Furthermore, let the closed neighbourhood $N_G[v]$ of $v$ be defined by $N_G[v]=N_G(v)\cup \{v\}$.
The maximum degree $\Delta(G)$ of $G$ is the maximum over all vertex degrees in $G$.
The minimum degree $\delta(G)$ of $G$ is the minimum over all vertex degrees in $G$.
For a set $X\subseteq V(G)$, let $G[X]$ be the subgraph induced by $X$.
We write $G\setminus X$ for $G[V(G)\setminus X]$.
If $v\in V(G)$ and $e\in E(G)$, we write $G-v$ for $G[V(G)\setminus \{v\}]$ and $G-e$ for the graph in which $e$ is deleted from $G$.
We call a set $I$ of vertices of $G$ an independent set of $G$, if no two vertices in $I$ are adjacent.
The maximum order of an independent set of $G$ is the independence number $\alpha(G)$.
A set $M$ of edges of $G$ is a matching, if no two edges in $M$ are adjacent.
The maximum size of a matching in $G$ is the matching number $\beta(G)$.
We say $G$ has a perfect matching, if $\beta(G)=\frac{1}{2}n(G)$.
We call a graph $G$ factor-critical, if for every vertex $v\in V(G)$, the graph $G-v$ has a perfect matching.
A bridge is an edge $e\in E(G)$ such that $G-e$ has more components than $G$.
We call a component $C$ of a graph odd (even), if $n(C)$ is odd (even).
Let $o(G)$ be the number of odd components of $G$.
For a graph $G$, we call a set $X$ of vertices of $G$ a clique, if all vertices in $X$ are pairwise adjacent.
The clique number $\omega(G)$ is the largest order of a clique in $G$.
The chromatic number $\chi(G)$ is the smallest integer $k$ such that
$V(G)$ has a partition into $k$ independent sets.

\bigskip
\noindent
Now we come to our main results.
\begin{thm}\label{mainthm}
If $G$ is a triangle-free graph with $\Delta(G)\leq 4$,
then
\begin{align*}
	\frac{7}{4}\alpha(G)+\beta(G) \geq n(G)
\end{align*}
with equality if and only if every component $C$ of $G$ has order $13$, $\alpha(C)=4$, and $\beta(C)=6$.
\end{thm}

\begin{thm}\label{mainthm2}
If $G$ is a triangle-free graph with $\Delta(G)\leq 4$,
then
\begin{align*}
	\alpha(G)+\frac{3}{2}\beta(G) \geq n(G)
\end{align*}
with equality if and only if every component $C$ of $G$ is either a single vertex,
or a cycle on $5$ vertices, or
has order $13$, $\alpha(C)=4$ and $\beta(C)=6$.
\end{thm}

\noindent
Theorem \ref{mainthm2} strengthens a result of Henning, L\"owenstein and Rautenbach by relaxing the maximum degree condition from $3$ to $4$.
Theorem \ref{mainthm} lead to some consequences for $\chi$-binding functions.
We say that a class $\mathcal{G}$ of graphs has a $\chi$-binding function if there is a function $f:\mathbb{N}_0\rightarrow \mathbb{N}_0$ such that 
$\chi(G)\leq f(\omega(G))$ for every $G\in \mathcal{G}$.
As mentioned in \cite{helora} the class of graphs that satisfy $\alpha(G)\leq 2$ has the $\chi$-binding function $\chi(G)=O \left( \frac{\omega(G)^2}{\log \omega(G)}\right)$.

Choudum et al. \cite{ch} give $\chi$-binding functions for graphs with independence number at most 2 and an excluded graph of order $5$.
Our main result implies the following.
\begin{coro}
If $G$ is a graph that contains neither $3K_1$ nor $K_1\cup K_5$ as an induced subgraph,
then 
\begin{align*}
	\chi(G)\leq \frac{7}{4}\omega(G).
\end{align*}
\end{coro}
\noindent
\textit{Proof:} The complement $\bar{G}$ of $G$ is triangle-free and satisfies $\Delta(\bar{G})\leq 4$.
Since there is no independent set of order at least 3, we obtain $\chi(G)=n(G)-\beta(\bar{G})$.
Since $\frac{7}{4}\alpha(\bar{G})+\beta(\bar{G})\geq n(\bar{G})=n(G)$ holds, 
it follows
\begin{align*}
	\chi(G)\leq \frac{7}{4}\alpha(\bar{G})=\frac{7}{4}\omega(G).
\end{align*}
$\Box$

\bigskip

\noindent
Our results are best possible.
Consider the graph $G_{13}$ on $13$ vertices $v_1,\ldots,v_{13}$.
Two vertices $v_i$ and $v_j$ such that $i<j$ are joined by an edge if and only if $j-i\in \{1,5,8,12\}$.
That means $G_{13}$ is a cycle with all $5$-chords.
See Figure 1 for an illustration.
Trivially $\beta(G_{13})=6$ and to convince ourselves that $\alpha(G_{13})=4$, we consider an independent set $I$ of size $\alpha(G_{13})$.
By symmetry we assume $v_1\in I$.
Let $G'=G_{13}\setminus N_{G_{13}}[v_1]$.
The remaining graph $G'$ is a cycle on $8$ vertices $v_3,v_4,v_5,v_{10},v_{11},v_{12},v_7,v_8$ with two 4-chords $v_3v_{11}$ and $v_4v_{12}$.
For contradiction we suppose $\alpha(G')\geq4$ and let $I'$ be an independent set of $G'$ of size $4$.
Since there is a spanning 8-cycle in $G'$, every second vertex of the cycle is in $I'$,
thus either $v_3$ and $v_{11}$ or $v_4$ and $v_{12}$ is in $I'$,
which is a contradiction.
It follows $\alpha(G') = 3$, because $\{v_3,v_5,v_7\}$ is independent in $G'$.
Since $\alpha(G_{13})= 1 + \alpha(G')$, by symmetry,
we conclude $\alpha(G_{13})=4$.
Thus $\frac{7}{4}\alpha(G_{13})+\beta(G_{13})=13=n(G_{13})$.

\begin{figure}[ht]
\begin{center}
\begin{tikzpicture}
\def\x{2}
\def\y{2.5}

\def\va{90}
\def\vb{117.7}
\def\vc{145.4}
\def\vd{173.1}
\def\ve{200.8}
\def\vf{228.5}
\def\vg{256.2}
\def\vh{283.9}
\def\vi{311.6}
\def\vj{339.3}
\def\vk{7}
\def\vl{35.7}
\def\vm{62.4}

\draw[thick] (0,0) +(\va:\x cm)-- +(\vb:\x cm)-- +(\vc:\x cm)-- +(\vd:\x cm)--
 +(\ve:\x cm)-- +(\vf:\x cm)-- +(\vg:\x cm)-- +(\vh:\x cm)-- +(\vi:\x cm)-- +(\vj:\x cm)
 -- +(\vk:\x cm)-- +(\vl:\x cm)-- +(\vm:\x cm) -- +(\va:\x cm);

\draw[thick] (0,0) +(\va:\x cm)-- +(\vf:\x cm)-- +(\vk:\x cm)-- +(\vc:\x cm)--
 +(\vh:\x cm)-- +(\vm:\x cm)-- +(\ve:\x cm)-- +(\vj:\x cm)-- +(\vb:\x cm)-- +(\vg:\x cm)
 -- +(\vl:\x cm)-- +(\vd:\x cm)-- +(\vi:\x cm) -- +(\va:\x cm);

\path[fill] (\va:\x) circle (0.5ex);
\path[fill] (\vb:\x) circle (0.5ex);
\path[fill] (\vc:\x) circle (0.5ex);
\path[fill] (\vd:\x) circle (0.5ex);
\path[fill] (\ve:\x) circle (0.5ex);
\path[fill] (\vf:\x) circle (0.5ex);
\path[fill] (\vg:\x) circle (0.5ex);
\path[fill] (\vh:\x) circle (0.5ex);
\path[fill] (\vi:\x) circle (0.5ex);
\path[fill] (\vj:\x) circle (0.5ex);
\path[fill] (\vk:\x) circle (0.5ex);
\path[fill] (\vl:\x) circle (0.5ex);
\path[fill] (\vm:\x) circle (0.5ex);

\node (1) at (\va:\y) {$v_1$};
\node (2) at (\vb:\y) {$v_2$};
\node (3) at (\vc:\y) {$v_3$};
\node (4) at (\vd:\y) {$v_4$};
\node (5) at (\ve:\y) {$v_5$};
\node (6) at (\vf:\y) {$v_6$};
\node (7) at (\vg:\y) {$v_7$};
\node (8) at (\vh:\y) {$v_8$};
\node (9) at (\vi:\y) {$v_9$};
\node (10) at (\vj:\y) {$v_{10}$};
\node (11) at (\vk:\y) {$v_{11}$};
\node (12) at (\vl:\y) {$v_{12}$};
\node (13) at (\vm:\y) {$v_{13}$};

\end{tikzpicture}
\end{center}
\caption{$G_{13}$}
\end{figure}

\noindent
Before we prove Theorem \ref{mainthm} and Theorem \ref{mainthm2} in Section 2 and 3,
we recall a few results for later use.
The first result due to Jones \cite{jo}.

\begin{thm}[Jones \cite{jo}]\label{indratio}
If $G$ is a triangle-free graph with $\Delta(G)\leq 4$,
then $$\alpha(G)\geq \frac{4}{13}n(G).$$
\end{thm}

\noindent
This result leads to a simple consequence.

\begin{coro}\label{perfectM}
Let $G$ be a triangle-free graph such that $\Delta(G)\leq 4$.
If $\beta(G)> \frac{6}{13}n(G)$, then 
\begin{align*}
	\frac{7}{4}\alpha(G)+\beta(G) > n(G) \quad {\rm and}\quad \alpha(G)+\frac{3}{2}\beta(G) > n(G).
\end{align*}
\end{coro}

\noindent
In addition we need Gallai and Edmonds' theorem on matchings in graphs.

\begin{thm}[Gallai \cite{ga} and Edmonds \cite{ed}]\label{gallai}
If $G$ is a graph, then there is a set $X\subseteq V(G)$ such that
\begin{align*}
	\beta(G)=\frac{1}{2}(n(G)+|X|-o(G\setminus X))
\end{align*}
and every odd component of $G\setminus X$ is factor-critical.
\end{thm}

\section{Proof of Theorem \ref{mainthm}}
\noindent
For contradiction, we suppose that $G$ is a counterexample such that the order of $G$ is as small as possible.
Because of the linearity of $\alpha(G)$ and $\beta(G)$ with respect to the components, $G$ is connected.
Note that every triangle-free graph $G$ has the following property:
\begin{align*}
\alpha(G)\geq \left\{
\begin{array}{rl}
	1, & \text{if } n(G)\geq 1,\\
	2, & \text{if } n(G)\geq 3,\\
	3, & \text{if } n(G)\geq 6, \text{ and}\\
	4, & \text{if } n(G)\geq 9.
\end{array}\right.
\end{align*}
This is a direct consequence of the Ramsey numbers $r(3,i)$ for $i\in \{1,2,3,4\}$ which are $r(3,1)=1$, $r(3,2)=3$, $r(3,3)=6$, and $r(3,4)=9$.
In addition there is no triangle-free factor-critical graph of order $3$.
By Theorem \ref{gallai}, there is some $X\subseteq V(G)$ such that
\begin{align}\label{beta1}
	\beta(G)=\frac{1}{2}(n(G)+|X|-o(G\setminus X))
\end{align}
and every component of $G\setminus X$ is factor-critical.
\bigskip

\noindent
\textbf{Claim 1:}
$\delta(G)\geq2$.
\bigskip

\noindent
\textit{Proof of Claim 1:} Trivially we obtain $\delta(G)\geq 1$.
Let us suppose for contradiction that there is a vertex $v\in V(G)$ of degree 1.
Let $u$ be the neighbour of $v$ and let $G'=G\setminus \{u,v\}$.
Since $\frac{7}{4}\alpha(G')+\beta(G')\geq n(G')=n(G)-2$, $\alpha(G)\geq 1+\alpha(G')$ and $\beta(G)\geq 1 + \beta(G')$,
we conclude
\begin{align*}
	\frac{7}{4}\alpha(G)+\beta(G)\geq n(G) \geq \frac{7}{4}\alpha(G')+\beta(G') + \frac{11}{4} > n(G).
\end{align*}
This is a contradiction and completes the proof. $\Box$

\bigskip

\noindent
\textbf{Claim 2:}
$n(G)\geq 14$.
\bigskip

\noindent
\textit{Proof of Claim 2:} We leave it to the reader to verify the case $n(G)\leq 5$.
If $X$ is empty, then either $G$ is factor-critical or has a perfect matching and this implies the claim by using the above Ramsey numbers.
Thus we assume $X\neq \emptyset$.
Because $G$ is connected, there is at least one edge joining every component of $G\setminus X$ to $X$.
Note that $\alpha(G)\geq o(G\setminus X)=n(G)+|X|-2\beta(G)$.
Since $\beta(G)< \frac{3}{10}n(G)+ \frac{7}{10}$ is equivalent to $\frac{7}{4}(n(G)+1-2\beta(G))+\beta(G)>n(G)$,
we assume
\begin{align}\label{beta}
	\beta(G)\geq \frac{3}{10}n(G)+ \frac{7}{10}.
\end{align}
By Claim $1$, there is no odd component of order $1$, if $|X|=1$.
Let us now suppose $6\leq n(G) \leq 8$. Now $\alpha(G)\geq3$ and (\ref{beta}) imply the claim.

Let us now suppose $9\leq n(G) \leq 12$.
This leads to $\alpha(G)\geq \max\{4,n(G)+1-2\beta(G)\}$.
If $n(G)\in\{9,10\}$, then $\alpha(G)\geq4$ and (\ref{beta}) imply the claim.

If $n(G)=11$, then $\alpha(G)\geq4$ and (\ref{beta}) imply $\frac{7}{4}\alpha(G)+\beta(G) \geq n(G)$.
To prove that this inequality is strict, 
we suppose for contradiction $\alpha(G)=4$ and $\beta(G)=4$.
This implies $|X|=1$ and $o(G\setminus X)=4$ and hence $n(G)\geq 21$, 
which is a contradiction.

If $n(G)=12$, then (\ref{beta}) implies $\beta(G)\geq 5$.
Since $\alpha(G)\geq 4$, it follows $\frac{7}{4}\alpha(G)+\beta(G) \geq n(G)$.
Again we suppose for contradiction $\alpha(G)=4$ and $\beta(G)=5$.
If $|X|=1$, then $o(G\setminus X)=3$.
This implies $n(G)\geq 16$, which is a contradiction.
If $|X|=2$, then $o(G\setminus X)=4$.
Thus there is an even component or an odd component of order at least $5$.
If $|X|\geq3$, then $o(G\setminus X)\geq 5$.
Both leads to $\alpha(G)\geq 5$, which is a contradiction. 

If $n(G)=13$, then (\ref{beta}) implies $\beta(G)\geq 5$.
If $\beta(G)=6$ or $\alpha(G)\geq 5$, then $G$ is not a counterexample.
Thus we suppose $\beta(G)=5$ and $\alpha(G)=4$.
If $|X|=1$, then $o(G\setminus X)=4$ and hence $n(G)\geq 21$.
If $|X|\geq 2$, then $o(G\setminus X)\geq 5$ and thus $\alpha(G)\geq 5$.
In both cases this is a contradiction, which completes the proof of the claim.
$\Box$
\bigskip

\noindent
\textbf{Claim 3:}
\textit{There is no bridge $e$ of $G$ such that $G - e$ has a factor-critical component.}
\bigskip

\noindent
\textit{Proof of Claim 3:}
For contradiction, we suppose there is a bridge $e=uv$ such that there is a factor-critical component $C$ in $G- e$.
Say $u\in V(C)$.
Let $G'=G[V(G)\setminus(V(C)\cup \{v\})]$.
Note that $G[V(C)\cup \{v\}]$ has a perfect matching and $\frac{7}{4}\alpha(G')+\beta(G')\geq n(G')$.
Since $\alpha(G)\geq \alpha (G')+\alpha(C)$
and $\beta(G) \geq \beta(G')+ \frac{1}{2}(n(C)+1)$, 
we obtain
\begin{align*}
	\frac{7}{4}\alpha(G)+\beta(G)
	&\geq \frac{7}{4}\left(\alpha(G')+ \alpha(C)\right)+\beta(G')+ \frac{1}{2}(n(C)+1)\\
	&\geq n(G')+ \frac{7}{4}\alpha(C)+ \frac{1}{2}n(C)+ \frac{1}{2}\\
	&= n(G)+ \frac{7}{4}\alpha(C)- \frac{1}{2}n(C)- \frac{1}{2}.
\end{align*}
Let $f(C)=\frac{7}{4}\alpha(C)- \frac{1}{2}n(C)- \frac{1}{2}$.
It remains to show that $f(C)>0$.
If $n(C)=1$, then $f(C)=\frac{3}{4}$.
If $n(C)=5$, then $f(C)\geq \frac{1}{2}$.
If $n(C)=7$, then $f(C)\geq \frac{5}{4}$.
If $n(C)=9$, then $f(C)\geq 2$.
If $n(C)=11$, then $f(C)\geq 1$.
If $n(C)=13$, then $f(C)\geq \frac{7}{4}$, 
because every triangle-free graph on $13$ vertices with at least one vertex of degree at most $3$ has independence number at least $5$.
To see this, let $v\in V(C)$ be such that $d_C(v)\leq 3$ and $\tilde{C}=C\setminus N_C[v]$.
We have $\alpha(C)\geq \alpha(\tilde{C})+1$.
Since $n(\tilde{C})\geq 9$, we obtain $\alpha(\tilde{C})\geq 4$ and the desired result follows.
If $n(C)\geq 15$, we conclude with Theorem \ref{indratio} 
\begin{align*}
	f(C)\geq \left(\frac{7}{13}- \frac{1}{2}\right)n(C) - \frac{1}{2}= \frac{1}{26}n(C)- \frac{1}{2}>0,
\end{align*}
which proves the claim.
$\Box$
\bigskip

\noindent
First we suppose that $X$ is empty.
By Corollary \ref{perfectM}, $G$ has no perfect matching.
This implies that $G$ is factor-critical and hence $\beta(G)=\frac{n(G)-1}{2}$.
By Corollary \ref{perfectM}, $\beta\leq \frac{6}{13}n(G)$ and
hence $\frac{n(G)-1}{2}\leq \frac{6}{13}n(G)$,
which is equivalent to $n(G)\leq 13$.
This is a contradiction to Claim $3$.

\noindent
Now we suppose $X$ is not empty.
Let $c_i$, and $c_{\geq i}$ be the number of odd components of $G\setminus X$ of order $i$ and at least $i$, respectively.
Let $n(c_{\geq i})$ be the number of vertices in odd components of $G\setminus X$ of order at least $i$.
Let $R$ be the set of vertices that are in no odd component of $G\setminus X$ and not in $X$.
By Claim 2, there are at least $2$ edges joining an odd componet of $G\setminus X$ to $X$.
By double-counting the edges between $X$ and the factor-critical components of $G\setminus X$,
it follows $4|X| \geq 2o(G\setminus X)$ and hence $|X|-o(G\setminus X) \geq - \frac{1}{2} o(G\setminus X)$.
We obtain
\begin{align}\label{beta_o}
	\beta(G)
	=\frac{1}{2}(n(G)+|X|-o(G\setminus X))
	\geq \frac{n(G)}{2}- \frac{1}{4}o(G\setminus X).
\end{align}
By Corollary \ref{perfectM}, we have
\begin{align}\label{beta2}
	\frac{6}{13}n(G)\geq \beta(G).
\end{align}
From (\ref{beta1}) and (\ref{beta2}) we conclude $n(G)-|X|\geq \frac{14}{13}n(G)-o(G\setminus X)$.
Since
\begin{align*}
	c_1+5c_5+7c_7+9c_9+n(c_{\geq11})+|R|=n(G)-|X|,
\end{align*}
this implies
\begin{align*}
	2c_1+6c_5+8c_7+10c_9+c_{\geq 11}+n(c_{\geq 11})+|R|\geq \frac{14}{13}n(G),
\end{align*}
which is equivalent to
\begin{align}\label{ineq}
	0\geq \frac{n(G)}{2}- \frac{13}{14}c_1- \frac{39}{14}c_5 - \frac{26}{7}c_7- \frac{65}{14}c_9 - \frac{13}{28}c_{\geq 11}- \frac{13}{28}n(c_{\geq 11})- \frac{13}{28}|R|.
\end{align}
Note that 
\begin{align}\label{indset}
	\alpha(G)\geq \alpha(R)+ \sum_{C:\ C{\rm\ odd\ component\ of\ }G\setminus X}\alpha(C).
\end{align}
We come to our final conclusion. 
Using Theorem \ref{indratio}, our remarks on the order of independent sets in small graphs, and (\ref{indset}),
we obtain
\begin{eqnarray*}
	 &&\frac{7}{4}\alpha(G)+\beta(G)\\
	 &\stackrel{(\ref{beta_o}),(\ref{indset})}{\geq}&
	  \frac{7}{4}\left(c_1+2c_5+3c_7+ 4c_9+\frac{4}{13}n(c_{\geq11})+ \frac{4}{13}|R|\right)+ \frac{n(G)}{2} - \frac{1}{4}o(G\setminus X)\\
	 &=&
	 \frac{n(G)}{2} + \frac{3}{2}c_1+ \frac{13}{4}c_5+5c_7+ \frac{27}{4}c_9- \frac{1}{4}c_{\geq11}+ \frac{7}{13}n(c_{\geq 11})+ \frac{7}{13}|R|\\
	  &\stackrel{(\ref{ineq})}{\geq} &
	 n(G) + \frac{4}{7}c_1+ \frac{13}{28}c_5+ \frac{9}{7}c_7 + \frac{59}{28}c_9 - \frac{5}{7}c_{\geq 11} + \frac{27}{364}n(c_{\geq 11})+ \frac{27}{364}|R|.
\end{eqnarray*}
Since $n(c_{11})\geq 11c_{11}$, it follows
\begin{align*}
	- \frac{5}{7}c_{\geq 11} + \frac{27}{364}n(c_{\geq 11}) \geq \left(-\frac{5}{7}+11 \cdot\frac{27}{364}\right)c_{\geq11}= \frac{37}{364}c_{\geq11}\geq 0.
\end{align*}
In addition, we have
\begin{align*}
	\frac{4}{7}c_1+ \frac{13}{28}c_5+ \frac{9}{7}c_7 + \frac{59}{28}c_9 +\frac{37}{364}c_{\geq11}+ \frac{27}{364}|R|>0.
\end{align*}
This implies the desiered inequality as well as the statement about the extremal graphs,
which completes the proof.
$\Box$
\bigskip

\section{Proof of Theorem \ref{mainthm2}}
\noindent
This proof is similar to the proof of Theorem \ref{mainthm},
thus we just sketch the important steps.
Again let $X\subseteq V(G)$ be again such that
\begin{align*}
	\beta(G)=\frac{1}{2}(n(G)+|X|-o(G\setminus X))
\end{align*}
and every component of $G\setminus X$ is factor-critical.
\bigskip

\noindent
\textbf{Claim 1:}
$\delta(G)\geq2$.
\bigskip

\noindent
\textit{Proof of Claim 1:} That is analogue to Claim 1 of the last section. $\Box$
\bigskip

\noindent
\textbf{Claim 2:}
$n(G)\geq 14$.
\bigskip

\noindent
\textit{Proof of Claim 2:}
Claim 2 can be proven with the same arguments as before.
The details are left to the reader.
$\Box$

\bigskip

\noindent
\textbf{Claim 3:}
\textit{There is no bridge $e$ of $G$ such that $G - e$ has a factor-critical component.}
\bigskip

\noindent
\textit{Proof of Claim 3:}
For contradiction, we suppose there is a bridge $e=uv$ such that there is a factor-critical component $C$ in $G- e$.
Say $u\in V(C)$.
Let $G'=G[V(G)\setminus(V(C)\cup \{v\})]$.
Note that $G[V(C)\cup \{v\}]$ has a perfect matching and $\alpha(G')+\frac{3}{2}\beta(G')\geq n(G')$.
Since $\alpha(G)\geq \alpha (G')+\alpha(C)$
and $\beta(G) \geq \beta(G')+ \frac{1}{2}(n(C)+1)$, 
we obtain
\begin{align*}
	\alpha(G)+\frac{3}{2}\beta(G)
	&\geq \alpha(G')+ \alpha(C)+\frac{3}{2}(\beta(G')+ \frac{1}{2}(n(C)+1))\\
	&\geq n(G')+ \alpha(C)+ \frac{3}{4}n(C)+ \frac{3}{4}\\
	&= n(G)+ \alpha(C)- \frac{1}{4}n(C)- \frac{1}{4}.
\end{align*}
Let $f(C)=\alpha(C)- \frac{1}{4}n(C)- \frac{1}{4}$.
It remains to show that $f(C)>0$.
If $n(C)=1$, then $f(C)=\frac{1}{2}$.
If $n(C)\geq 5$, we conclude by Theorem \ref{indratio} 
\begin{align*}
	f(C)\geq \left(\frac{4}{13}- \frac{1}{4}\right)n(C) - \frac{1}{4}= \frac{3}{52}n(C)- \frac{1}{4}>0,
\end{align*}
which proves the claim.
$\Box$
\bigskip

\noindent
First we suppose that $X$ is empty.
There are two possibilities.
By Corollary \ref{perfectM}, $G$ has no perfect matching.
This implies, that $G$ is factor-critical and hence $\beta(G)=\frac{n(G)-1}{2}$.
By Corollary \ref{perfectM}, $\beta\leq \frac{6}{13}n(G)$ and
hence $\frac{n(G)-1}{2}\leq \frac{6}{13}n(G)$,
which is equivalent to $n(G)\leq 13$.
This is a contradiction to Claim $3$.

\noindent
Now we suppose $X$ is not empty.
We use the same terminology as in the proof of Theorem $1$.
The inqualities (\ref{beta_o}), (\ref{ineq}), and (\ref{indset}) are still true.
We restate (\ref{ineq}) in an equivalent form
\begin{align}\label{ineq2}
	0\geq \frac{n(G)}{4}- \frac{13}{28}c_1- \frac{39}{28}c_5 - \frac{13}{7}c_7- \frac{65}{28}c_9 - \frac{13}{56}c_{\geq 11}- \frac{13}{56}n(c_{\geq 11})- \frac{13}{56}|R|.
\end{align}
We conlude
\begin{eqnarray*}
	&&\alpha(G)+\frac{3}{2}\beta(G)\\
	&\stackrel{(\ref{beta_o}),(\ref{indset})}{\geq} & 
	c_1+2c_5+3c_7+ 4c_9+\frac{4}{13}n(c_{\geq11})+ \frac{4}{13}|R|+ \frac{3}{2}\left(\frac{n(G)}{2} - \frac{1}{4}o(G\setminus X)\right)\\
	&=& 
	\frac{3n(G)}{4} + \frac{5}{8}c_1+ \frac{13}{8}c_5+ \frac{21}{8}c_7+ \frac{29}{8}c_9- \frac{3}{8}c_{\geq11}+ \frac{4}{13}n(c_{\geq 11})+ \frac{4}{13}|R|\\
	& \stackrel{(\ref{ineq2})}{\geq}& 
	n(G) + \frac{9}{56}c_1+ \frac{13}{56}c_5+ \frac{43}{56}c_7 + \frac{73}{56}c_9 - \frac{17}{28}c_{\geq 11} + \frac{55}{728}n(c_{\geq 11})+ \frac{55}{728}|R|.
\end{eqnarray*}
Since $n(c_{11})\geq 11c_{11}$, we have
\begin{align*}
	- \frac{17}{28}c_{\geq 11} + \frac{55}{728}n(c_{\geq 11}) \geq \left(-\frac{17}{28}+11 \cdot\frac{55}{728}\right)c_{\geq11}= \frac{163}{728}c_{\geq11}\geq 0.
\end{align*}
In addition, we have
\begin{align*}
	\frac{9}{56}c_1+ \frac{13}{56}c_5+ \frac{43}{56}c_7 + \frac{73}{56}c_9 + \frac{163}{728}c_{\geq 11}+ \frac{55}{728}|R|>0.
\end{align*}
This implies the desiered inequality as well as the statement about the extremal graphs,
which completes the proof.
$\Box$

\bigskip


\end{document}